\documentclass[10pt,a4paper,draft]{article}
\date{\today}

\usepackage{latexsym,amssymb,amsmath}
\usepackage[latin1]{inputenc}
\usepackage[T1]{fontenc}
\usepackage[english]{babel}

\def\RR{\mbox{\it I\hskip -0.177em R}}
\def\PP{\mbox{\it I\hskip -0.177em P}}
\def\NN{\mbox{\it I\hskip -0.177em N}}

\def\EE{\mbox{\it I\hskip -0.177em E}}
\def\II{\mbox{\rm \Large 1\hskip -0.353em 1}}
\newcommand{\sRR}{{\hbox{$\scriptstyle{I}$\kern-.25em\hbox{$\scriptstyle{R}$}}}}
\newcommand{\sPP}{{\hbox{$\scriptstyle{I}$\kern-.25em\hbox{$\scriptstyle{P}$}}}}
\newcommand{\sNN}{{\hbox{$\scriptstyle{I}$\kern-.25em\hbox{$\scriptstyle{N}$}}}}
\newcommand{\sZZ}{{\hbox{$\scriptstyle{Z}$\kern-.3em\hbox{$\scriptstyle{Z}$}}}}
\newcommand{\sEE}{{\hbox{$\scriptstyle{I}$\kern-.25em\hbox{$\scriptstyle{E}$}}}}
\newcommand{\sII}{{\hbox{$\scriptstyle{I}$\kern-.25em\hbox{$\scriptstyle{I}$}}}}

\RequirePackage{theorem}

\newtheorem{thm}{\noindent Theorem}[section]

\newtheorem{coro}{\noindent Corollary}[section]

\newtheorem{conj}{\noindent Conjecture}[section]
\newtheorem{rem}{\noindent Remark}[section]

\theoremheaderfont{\normalfont\bfseries} \theorembodyfont{\slshape}
{\theorembodyfont{\rmfamily} \newtheorem{ex}{\noindent
Example}[section]} {\theorembodyfont{\rmfamily}
}
{\theorembodyfont{\rmfamily} } {\theorembodyfont{\rmfamily}
}

\newenvironment{proof}{
\noindent {\bf Proof} \noindent} {\hfill $\Box$\vskip 5mm}

\begin{document}

\title{Limiting search cost distribution for the move-to-front\\ rule
with random request probabilities} 


\author{Javiera Barrera$^{\text{a}}$\footnote{Corresponding author} \and Thierry Huillet$^{\text{b}}$ \and Christian Paroissin$^{\text{c}}$}

\date{}

\maketitle

\begin{center}
\begin{trivlist}
\item {(a)} MAP5 - UMR 8145 CNRS, Universit\'e Paris 5-Ren\'e Descartes, 45 rue
des Saints-P\`eres, 75270 PARIS cedex 06, FRANCE.

\item{(b)} LPTM - UMR 8089 CNRS, Universit\'e de Cergy-Pontoise, 2
avenue Adolphe-Chauvin, 95302 CERGY-PONTOISE, FRANCE.

\item{(c)} LMA - UMR 5142 CNRS, Universit\'e de Pau et des Pays de l'Adour,
Avenue de l'Universit\'e, BP 1155, 64013 PAU cedex, FRANCE.
\end{trivlist}
\end{center}

\begin{abstract}
Consider a list of $n$ files whose popularities are random. These
files are updated according to the move-to-front rule and we
consider the induced Markov chain at equilibrium. We give the exact 
limiting distribution of the search-cost per item as $n$ tends to infinity.
Some examples are supplied. 
\end{abstract}



 \footnotetext[1]{E-mail Addresses: {\tt
jbarrera@dim.uchile.cl} (J. Barrera), {\tt
Thierry.Huillet@ptm.u-cergy.fr} (T. Huillet) and {\tt
cparoiss@univ-pau.fr} (C. Paroissin).}

\newpage
\section{Introduction and model}
\label{sec:intro}

Consider a list of $n$ files which is updated as follows: at each
unit of discrete time, a file is requested independently of the
previous requests and is moved to the front of the list. This
heuristic is called the move-to-front rule and was first introduced
by \cite{tsetlin} and \cite{mccabe} to sort files. Such strategy
is used when the request probabilities are unknown,
otherwise we would list the files in order to have decreasing
request probabilities. The move-to-front rule induces a Markov chain
over the permutations of $n$ elements which has a unique stationary
distribution, (see \cite{donnelly} and reference to the work of
Hendricks, Dies and Letac therein). This distribution turns out to
be the size-biased permutation of the request probabilities.
 \\[1ex]
Here, we consider that these request probabilities are themselves
random, as in a Bayesian analysis. Let $\omega=(\omega_i)_{i \in
\sNN^*}$ be a sequence of iid positive random variables. The
Laplace transform of a weight will be denoted by $\phi$ and its
expectation by $\mu$. For any $i\in \NN^*$, $\omega_i$ represents 
the weight of the file $i$. We can construct request probabilities 
${\bf p}=(p_1,..,p_n)$ as follows:
\begin{equation*}
\forall i\in \{1,...,n\}\;, \qquad p_i=\frac{\omega_i}{W_n}\quad
{\mbox{where}} \quad W_n=\sum_{i=1}^n \omega_i \;.
\end{equation*}
Such random vector ${\bf p}$ is called a random discrete
distribution \cite{kingman}.
\\[1ex]
Let us denote by $S_n$ the search cost of an item (i.e. the position
in the list of the requested item) when the underlying Markov chain
is in steady state (the first position will be 0). For this model, 
\cite{barreraparoissin}
obtained exact and asymptotic formulae for the Laplace transform of
$S_n$ (some results were also extended to the case of independent
random weights). In particular, they found the limit of the
expectation and the variance of $S_n$. Moreover, in the case of
i.i.d. gamma weights, \cite{barrerahuilletparoissin}, obtained the
exact and asymptotic distribution of $S_n$, using an exact
representation of the size-biased permutation arising from Dirichlet
partitions. Note that \cite{fillmtf} found the limiting
distribution of $S_n$ when weights are deterministic but
non-identical, in some cases (uniform, Zipf's law, generalized
Zipf's law, power law and geometric).
\\[2ex]
In section \ref{sec:limit}, we shall give a general formula
for the density of the limiting search cost distribution $S$, provided
that the expected weight is finite. Then we derive the moment function 
and the cumulative distribution function of $S$. We also discuss the 
relationship between the move-to-front rule and the least-recently-used
strategy. In section \ref{sec:ex} we study
some examples for which computations can be done explicitly: both
continuous and discrete distributions are considered.

\section{Limiting search cost distribution}\label{sec:limit}

The early analysis of the heuristic move-to-front
focused on the expected search cost, see \cite{mccabe},
\cite{kingman} and \cite{FGT}, for instance. Later,
researchers paid much attention to the (transient and stationary)
distribution of the search cost (\cite{FillHolst}). Some of them
investigated the limiting behavior as the number $n$ of items tends
to infinity (see \cite{fillmtf}). In a more recent article,
\cite{barreraparoissin} obtained an integral representation of the
Laplace transform of $S_n$ in the Bayesian model described in the
introduction. Their main theorem is the following:
\begin{thm}\label{thm:mtf2004}
For a sequence $\omega$ of iid positive random variables,
\begin{equation*}
\forall s \geqslant 0\,,\qquad \phi_{S_n}(s)= n  \int_0^{\infty}\int_t^{\infty}
\phi''(r)\left[ \phi(r)+e^{-s}\left(\phi(r-t)-\phi(r)\right)\right]^{n-1}
\,dr\,dt \;.
\end{equation*}
\end{thm}

In the same article the integral representation for the two first
moments of $S_n$ were derived. Moreover, they obtained a point-wise asymptotic
equivalent for the Laplace transform of $S_n$ and the limit of the
first two moments of $S_n/n$ when the number $n$ of items tends to
infinity. From theorem \ref{thm:mtf2004}, we can obtain the following 
closed-form expression for the density function of the limiting 
distribution of $Sn/n$:
\begin{thm}\label{thm:limdib}
For a sequence $\omega$ of iid positive random weights with finite 
expectation $\mu$,
\begin{equation*}
\frac{S_n}{n} \xrightarrow[n \rightarrow
  \infty]{d} S \;,
\end{equation*}
where $S$ is a continuous random variable with the following density
function $f_{S}$:
\begin{equation}\label{eqn:limdib}
f_{S} (x) = -\frac{1}{\mu}
\frac{\phi''\left(\phi^{-1}(1-x)\right)}{\phi'\left(\phi^{-1}(1-x)\right)}\II_{[0,1-p_0]}(x)
\;,
\end{equation}
where $p_0=\PP(\omega_i=0)$ and $\phi^{-1}$ is the inverse
function of $\phi$.
\end{thm}

\begin{rem}\label{rm_p_0}
The quantity $p_0$ can be interpreted as follows: $p_0$ is the
probability that an item is never requested. At stationarity, one
 expects that any such item will be at the bottom of the list: $np_0$
is the mean number of unrequested items. So it is not surprising
that the support of $S$ is not the entire unit interval. Note that
if the distribution of the weight is continuous, then $p_0=0$.
\end{rem}

\begin{proof}
We have to prove that $S_n/n$ converges in distribution, as
$n$ tends to infinity, to a certain random variable that will be
denote by $S$. First, observe that:
\begin{equation*}
\forall s \geqslant 0\,, \quad \phi_{S_n/n}(s) = \phi_{S_n}\left(\tfrac{s}{n}
\right) \;.
\end{equation*}
So we are now interested in the limit of $\phi_{S_n}(s/n)$.
\\[1ex]
For any reals $a$ and $b$ such that $0 \leqslant a \leqslant b \leqslant 
\infty$, let:
\begin{equation*}
I_n(a,b)=\int^b_a\phi''(r)\left[\phi(r)+e^{-s/n}(\phi(r-t)-\phi(r))
\right]^{n-1}\, dr\,.
\end{equation*}
If $b=\infty$, then we will omit this parameter, i.e. $I_n(a)=I_n(a,\infty)$. Using these notations, theorem \ref{thm:mtf2004} gives:
\begin{equation}\label{S_n}
\phi_{S_n}\left(\tfrac{s}{n}\right) = n\int_0^{\infty} I_n(t)\,dt \;.
\end{equation}
We now decompose $I_n(t)$ into two parts: $I_n(t) = I_n(t,t+\varepsilon)+
I_n(t+\varepsilon)$. We will prove that $nI_n(t+\varepsilon,\infty)$ 
tends to $0$ when $n$ tends to infinity:
\begin{eqnarray*}
n I_n(t+\varepsilon,\infty) & = & n\int_{t +\varepsilon}^{\infty}
\phi''(r)\left[e^{-s/n}(\phi(r-t)+ (1-e^{-s/n})\phi(r))\right]^{n-1} dr \;,\\
& \leqslant & n\int_{t +\varepsilon}^{\infty} \phi''(r) \phi(r-t)^{n-1}\,
dr\
,\\
& \leqslant & -n\phi(\varepsilon)^{n-1}\phi'(t+\varepsilon)\,,
\end{eqnarray*}
since $\phi$ is decreasing. Then $\lim_{n\rightarrow \infty}n 
I_n(t+\varepsilon,\infty)=0$, for all $\varepsilon>0$.
\\[1ex]
Now we will estimate $I_n(t,t+\varepsilon)$. Let $h_{n}(r,t) =
\phi(r)+e^{-s/n}(\phi(r-t)-\phi(r))$. For a fixed value of $t$, the function
$h_{n}(\cdot,t)$ behaves as $\phi$. In particular $\frac{\partial
h_{n}}{\partial r}$ is an increasing function for $r \in
[t,t+\varepsilon]$. Then we obtain the following bounds:
\begin{equation*}
\frac{\partial h_{n}}{\partial r}(t,t) \leqslant \frac{\partial
h_{n}}{\partial r}(r,t) \leqslant
\frac{\partial h_{n}}{\partial r}(t+\varepsilon,t)\,,
\end{equation*}
and:
\begin{equation*}
\phi''(t+\varepsilon) \leqslant \phi''(r) \leqslant \phi''(t)\,.
\end{equation*}

Hence, we can bound $I_n(t,t+\varepsilon)$ by:
\begin{eqnarray*}
I_n(t,t+\varepsilon) & = & \int_t^{t+\varepsilon} \phi''(r) \left(h_{n}(r,t)\right)^{n-1}\frac{\partial h_{n}}{\partial r}(r,t)\frac{\partial h_{n}}{\partial r}(r,t)^{-1} \, dr\\
& \leqslant &\phi''(t)\frac{\partial h_{n}}{\partial r}(t,t)^{-1}\int_t^{t+\varepsilon}  \left(h_{n}(r,t)\right)^{n-1}\frac{\partial h_{n}}{\partial r}(r,t) \, dr\\
& \leqslant &\phi''(t)\frac{\partial h_{n}}{\partial r}(t,t)^{-1}\frac{1}{n}
\left[\left(h_{n}(t+\varepsilon,t)\right)^n-\left(h_{n}(t,t)\right)^n\right]
\ .
\end{eqnarray*}
Proceeding similarly, we can find a lower bound:
\begin{equation*}
I_n(t,t+\varepsilon) \geqslant
\phi''(t+\varepsilon)\frac{\partial h_{n}}{\partial r}(t+\varepsilon,t)^{-1}\frac{1}{n}
\left[\left(h_{n}(t+\varepsilon,t)\right)^n-\left(h_{n}(t,t)\right)^n\right]\,.
\end{equation*}
Then, for any $\varepsilon>0$, one can prove the following limits hold:
\begin{eqnarray*}
\lim_{n \rightarrow \infty}\left(h_{n}(t+\varepsilon,t)\right)^n &
=  & 0\ ,\\
\lim_{n \rightarrow \infty}\left(h_{n}(t,t)\right)^n & =
&\exp\left[ -s(1-\phi(t)) \right]\ ,\\
\lim_{n \rightarrow \infty}\frac{\partial h_{n}}{\partial r}(t,t) & = &
\phi'(\varepsilon) \ ,\\
 \lim_{n \rightarrow
\infty}\frac{\partial h_{n}}{\partial r}(t+\varepsilon,t) & = & \phi'(0) \ .
\end{eqnarray*}
 Replacing these limits in the equations above, we have computed upper and lower bounds of $I_n (t,t+\varepsilon)$. In other words, if the limit of $nI_n (t,t+\varepsilon)$ exists, then it is bounded by:
\begin{equation*}
-\frac{\phi''(t+\varepsilon)}{\phi'(0)}
\exp\left( -(1-\phi(t))s \right) \leqslant
\lim_{n \rightarrow \infty} n I_n(t,t+\varepsilon) \leqslant
-\frac{\phi''(t)}{\phi'(\varepsilon)} \exp\left( -(1-\phi(t))s
\right)  \;.
\end{equation*}
This is true for any $\varepsilon>0$; then letting $\varepsilon$
tends to $0$, we have:
\begin{equation*}
\lim_{n \rightarrow \infty} n I_n(t)=
\frac{\phi''(t)}{\mu}\exp\left( -(1-\phi(t))s \right)\;.
\end{equation*}
Replacing this limit in equation (\ref{S_n}) we obtain
\begin{equation}\label{LaplaceT}
\lim_{n \rightarrow \infty} \phi_{S_n/n}(s)= \frac{1}{\mu} \int_0^{\infty} \phi''(t)  e^{-(1-\phi(t))s}\,dt \;,
\end{equation}
which will be denoted by $\phi_S(s)$. Although this limit a priori is not
necessarily the Laplace transform of a random variable, according to
the Continuity theorem (page 431 Ch. XIII in \cite{feller}), one has to check that $\lim_{s\rightarrow
  0}\phi_S(s)=1$, which can be proved by using the dominated convergence 
theorem.
\\[1ex]
A suitable change of variable $y=1-\phi(r)$ in equation (\ref{LaplaceT}) 
gives:
\begin{equation*}
\phi_{S}(s) = -\frac{1}{\mu} \int_0^{1-p_0}
\frac{\phi''\left(\phi^{-1}(1-y)\right)}{\phi'\left(\phi^{-1}(1-y)\right)}
e^{-ys}\, dr ,
\end{equation*}
where for the integral limits we used the property that
$\phi(\infty)= p_0$ (see \cite{feller} remark in theorem 1(a) page
439 Ch. XIII). Therefore, we have that:
\begin{equation*}
f_{S} (y) = -\frac{1}{\mu}
\frac{\phi''\left(\phi^{-1}(1-y)\right)}{\phi'\left(\phi^{-1}(1-y)\right)}\II_{[0,1-p_0]}(y)
\end{equation*}
is the probability density of $S$.
\end{proof}

As a corollary to this theorem, we can compute the $q$-th moment and the cumulative distribution function (c.d.f.) of $S$:
\begin{coro}
For any $q \in \RR$
\begin{equation*}
E[S^q] = \frac{1}{\mu} \int_{0}^{\infty} (1-\phi(t))^q
\phi''(t)\,dt\ ,
\end{equation*}
and, for any $x \in [0,1]$,
\begin{equation*}
\PP (S\leqslant x) = \left( \frac{1}{\mu} \int_{0}^{\phi^{-1}(1-x)} \phi''(t)\,dt \right)\II_{[0,1-p_0]}(x)+\II_{(1-p_0,1]}(x) \;.
\end{equation*}
\end{coro}

One could be interested in the cumulative distribution function of $S$ (or more precisely in the survival function), since the move-to-front rule is related to the least-recently-used strategy (see \cite{FGT} for instance). Indeed, many operating systems or softwares use a memory (also called cache) that could be quickly addressed (think of a web browser, for instance). Hence, one needs to define a strategy to organize it. Let us consider that the cache is made of $k$ files. The least-recently-used strategy is the following: at each unit of discrete time, a file is requested and is moved in front of the cache; if the file was not just previously in the cache, then the last file is deleted from the cache and all other files are shifted by one position to the right; if the file was  just previously in the cache, then the file is moved exactly as in the move-to-front rule. So, the move-to-front rule can be viewed as a special case of the least-recently-used strategy for which the length of the cache is equal to the number of files ($k=n$). An important question arises: what is the probability that the requested file is not in the cache? The probability of this event is called the page default; we will denote it by $\pi_k$ in the sequel. Because of the link between the move-to-front rule and the least-recently-used strategy (as underlined above), we clearly have that $\pi_k = \PP(S_n \geqslant k)$. So, if we assume that the cache length is proportional to the number of files, say $k=\alpha n$ with $\alpha \in [0,1]$ fixed, for a large collection of files, the following approximation holds:
\begin{equation*}
\pi_{\alpha n} \simeq \frac{1}{\mu} \int_{0}^{\phi^{-1}(1-\alpha)} \phi''(t)\,dt \;
\end{equation*}
if $\alpha<p_0$ and $\pi_{\alpha n} \simeq 1$ otherwise.

\section{Examples}
\label{sec:ex}

In this section, we study some examples for which we are able to do
explicitly all computations. We will consider both continuous and discrete
distribution for the random weights.

\begin{ex}
Suppose that the weights have the Dirac distribution at point mass
$1$ (in other words, weights are deterministic and are equally
requested). Then $\phi(r)=e^{-r}$, the expectation $\mu = 1$ and $
p_0=0$, we deduce that:
\begin{equation*}
f_{S_1}(x) = \II_{[0,1]}(x)\;.
\end{equation*}
Thus, $S_1$ has the uniform distribution over $[0,1]$: this result was
already proved in (theorem~4.2, p.~198 of \cite{fillmtf}). The
$k$-th moment (with $k \in \RR_+$) and the c.d.f. of $S_1$ is:
\begin{equation*}
\EE[S^k_1]= \frac{1}{k+1} 
\qquad {\mbox{and}} \qquad \forall x \in [0,1]\,, \quad
F_{S_1}(x) = \PP(S_1\leqslant x) = x  \;.
\end{equation*}
\end{ex}

\begin{ex}
Suppose that the weights have the Gamma distribution with parameter
$\alpha>0$. In this example, the random vector $(p_1, \ldots, p_n)$
has the symmetric Dirichlet distribution $D_n(\alpha)$ (see \cite{wilks} 
or \cite{kbj}). In such a 
case, $p_0=0$, $\mu = \alpha$ and $\phi(r)= (1+r)^{-\alpha}$.
Computations give:
\begin{equation*}
f_{S_2}(x) = \left(1+\frac{1}{\alpha}\right) (1-x)^{1/\alpha}
\II_{[0,1]}(x)\;,
\end{equation*}
which is the density function of the Beta distribution with
parameters $(1,1+1/\alpha)$. Note that this result has already been
proved by \cite{barrerahuilletparoissin} with a specific technique
using properties of Dirichlet distribution (in this case we were able not
only to find the limiting search cost distribution but also the transient
search cost distribution for any finite $n$). The $k$-th moment
(with $k \in \RR_+$) of $S_2$ is:
\begin{equation*}
\EE[S_2^k]=
\frac{\Gamma(k+1)\Gamma(2+\tfrac{1}{\alpha})}{\Gamma(2+k+\tfrac{1}{\alpha})}\;.
\end{equation*}
In particular, we have $\EE[S_2]=\tfrac{\alpha}{2\alpha+1}$ and
${\mbox{Var}}[S_2]= \tfrac{(\alpha+1)\alpha^2}{(3\alpha+1)(2\alpha+1)^2}$.
One can also compute the c.d.f of $S_2$ and, for any $x \in [0,1]$, we get:
\begin{equation*}
F_{S_2}(x) = \PP(S_2\leqslant x)= 1 - (1-x)^{1+1/\alpha}\;.
\end{equation*}
We can easily deduce that, for any $x \in [0,1]$, $\bar{F}_{S_2}(x) \leqslant 
\bar{F}_{S_1}(x)$, where $\bar{F}_{S_1}(\cdot) = 1 - F_{S_1}(\cdot)$. So we 
have $S_2 \preceq_{st} S_1$ (where $\preceq_{st}$ denotes the usual stochastic 
ordering; see \cite{shakedshanthikumar} or \cite{stoyan}, for instance).
\end{ex}

\begin{ex}
Suppose that the weights have the Geometric distribution on $\NN$
with parameter $p \in (0,1)$. In such case, $p_0=p$, $\mu = (1-p)/p$
and $\phi(r)=p/(1-(1-p)e^{-r})$. Elementary computations give:
\begin{equation*}
f_{S_3}(x)=\frac{2(1-x)-p}{1-p}  \II_{[0,1-p]}(x) \;.
\end{equation*}
The $k$-th moment (with $k \in \RR_+$) of $S_3$ is:
\begin{equation*}
\EE[S_3^k]= \frac{(2+pk)(1-p)^k}{(k+1)(k+2)} \;.
\end{equation*}
In particular, we have $\EE[S_3]=\tfrac{(2+p)(1-p)}{6}$ and
${\mbox{Var}}[S_3]= \tfrac{(1-p)^2(2+2p-p^2)}{36}$. One can also compute 
the c.d.f of $S_3$ and, for any $x \in [0,1]$, get:
\begin{equation*}
F_{S_3}(x) = \PP (S_3\leqslant x)= \frac{x(2-p-x)}{1-p}\II_{[0,1-p]}(x)+\II_{(1-p,1]}(x) \;.
\end{equation*}
Hence, from the above expression, one can check that $S_4 \preceq_{st} S_1$.
\end{ex}

\begin{ex}
Suppose that the weights have the Poisson distribution with
parameter $\lambda$. In such case, $p_0=e^{-\lambda}$, $\mu =
\lambda$ and $\phi(r)=\exp\left(\lambda e^{-r}-1\right)$. Simple
computations give:
\begin{equation*}
f_{S_4}(x)=\frac{\ln(1-x)+\lambda +
1}{\lambda}\II_{[0,1-e^{-\lambda}]}(x) \;.
\end{equation*}
Using formula~1.6.5.3 of \cite{prudnikovetal} (page~244), one can
compute the $k$-th moment (with $k \in \NN$) of $S_4$:
\begin{equation*}
\EE[S_4^k]= \frac{1}{\lambda(k+1)}\left[ \lambda + (1-e^{-\lambda})^{k+1} - \sum_{i=1}^{k+1}\frac{(1-e^{-\lambda})^i}{i} \right] \;.
\end{equation*}
In particular, we have $\EE[S_4]= \tfrac{1}{2} - 
\tfrac{1-e^{-2\lambda}}{4\lambda}$. One can also compute 
the c.d.f of $S_4$ and, for any $x \in [0,1]$, we get:
\begin{equation*}
F_{S_4}(x) = \PP (S_4\leqslant x)= (x-\frac{1}{\lambda}(1-x)\ln(1-x))\II_{[0,1-e^{-\lambda}]}(x)+\II_{(1-e^{-\lambda},1]}(x) \;.
\end{equation*}
Thus, from the expression above, one can deduce that $S_4 \preceq_{st} S_1$.
\end{ex}

From the study of these four examples, one can observe that both $S_2$, $S_3$ and $S_4$ are stochastically smaller than $S_1$. Hence, the following conjecture looks appealing:

\begin{conj}
Let $S$ be the limiting distribution of the search cost associated to a sequence $\omega$ of iid positive random variables. Then, $S \preceq_{st} S_1$ where $S_1$ is a random distribution having the uniform distribution on the unit interval.
\end{conj}

This conjecture is compatible with some remarks in \cite{barreraparoissin}, more precisely with proposition~3.1 therein. Indeed, if the conjecture is right, then as a consequence we have $\EE[S] \leqslant \EE[S_1] = \tfrac{1}{2}$. And this is precisely what is stated in  proposition~3.1. This conjecture can be interpreted as follows: the case with Dirac weights corresponds to the worst case. Despite our conjecture seems to be true, its proof seems to be difficult.


\vskip 0.5cm
\paragraph{Acknowledgment}{The authors thank support from FONDAP-CONICYT in Applied Mathematics and Millenium Nucleus in Information and Randomness ICM P01-005. JB wishes to thank CONICYT National Postgraduate Fellowship Program which supports her PhD.

\end{document}